\documentstyle[11pt, sec, psfig, diagram]{article}

\newtheorem{theorem}{Theorem}[section]
\newtheorem{lemma}[theorem]{Lemma}

\newtheorem{remark}[theorem]{Remark}

\newfont{\Bbb}{msbm10 scaled \magstep1}

\newcommand{\bQ}{{\hbox{\Bbb Q}}}
\newcommand\bR{\hbox{\Bbb R}}

\newcommand\bZ{\hbox{\Bbb Z}}

\newfont{\es}{eusm10 scaled \magstep1}
\newfont{\ses}{eufm8 scaled \magstep1}
\newfont{\gt}{eufb10 scaled \magstep1}
\newfont{\sg}{eufb8 scaled \magstep1}
\newfont{\goth}{eufb10 scaled \magstep2}

\newcommand{\p}{{\cal P}}

\def\ra{\rightarrow}

\def\be{\begin{equation}}
\def\ee{\end{equation}}

\newcommand{\LL}{\left(\left( }
\newcommand{\RR}{\right)\right)}

\newcommand{\modu}{\hbox{\gt M}}

\newcommand{\ve}{{\varepsilon}}

\begin{document}

\title{Lattice Points, Dedekind-Rademacher Sums and a Conjecture of Kronheimer and Mrowka}

\author{Liviu I.Nicolaescu\thanks{{\bf Current address}: Dept.of Math., McMaster University, Hamilton, Ontario,  L8S 4K1, Canada; nicolaes@icarus.math.mcmaster.ca}}

\date{Version 2, January 1998}

\maketitle

\begin{abstract} 
We express the number of lattice points inside a certain simplex with vertices in ${\bQ}^3$ or ${\bQ}^4$ in terms of Dedekind-Rademacher sums.  As an application  we prove   a conjecture  of Kronheimer and Mrowka  in the special case of Brieskorn homology spheres $\Sigma(a_1, \cdots, a_n)$, $n\leq 4$.  This conjecture relates the Euler characteristic of the Seiberg-Witten-Floer homology to the Casson invariant. \end{abstract}

\addcontentsline{toc}{section}{Introduction}

\setcounter{section}{0}
\begin{center}
{\bf Introduction}
\end{center}

\bigskip

Since the very beginning it was apparent that the Seiberg-Witten analogue of the instanton Floer homology of a homology 3-sphere is no longer a topological invariant since it can vary with the metric.

Recently, M. Marcolli  explained in \cite{Mar} the metric dependence of the Euler characteristic of the SWF ($=$ Seiberg-Witten-Floer)  homology.  More precisely, if $g_i$  ($i=0,1$) are two generic Riemann metrics on a homology $3$-sphere $N$ and $\chi_{SW}(N,g_i)$ is the Euler characteristic of the SWF homology of $(N,g_i)$  the results of  \cite{Mar} imply  that
\[
\chi_{SW}(N,g_1)-\chi_{SW}(N,g_0)= \frac{1}{8}{\bf F}(g_1)-\frac{1}{8}{\bf F}(g_0)
\]
where
\[
{\bf  F}(g)=  4\eta_{dir}(g)+\eta_{sign}(g)
\]
$\eta_{dir}(g)$ denotes the eta invariant of the Dirac operator of $(N,g)$ while $\eta_{sign}(g)$ denotes the eta invariant of the odd signature operator on $(N,g)$. In particular, the above equality shows that  the quantity
\[
\alpha(N) =-\chi_{SW}(N,g) +\frac{1}{8}{\bf F}(g)
\]
is independent of $g$  and thus is  a topological quantity.

In 1996, Kronheimer and Mrowka conjectured (see \cite{KM}) that this quantity coincides (up to a sign) with the Casson invariant of  $N$.   W. Chen (\cite{C1},\cite{C2}) has  established many interesting properties of this invariant. 

The topological goal of this paper is to provide a proof of this conjecture,  in the special case of Brieskorn homology sphere $\Sigma(a_1,\cdots a_n)$, $n\leq 4$. Our proof is arithmetic in nature and is based on a count of lattice points inside a certain $n$-simplex   in  ${\bR}^n$ with vertices in ${\bQ}^n$, $n=3,4$.

More precisely, recent work (\cite{MOY}, \cite{N}) shows that for a certain natural metric $g_0$ on $\Sigma(a,b,c)$  (realizing the Thurston geometry of this Seifert manifold) we have
\[
-\chi_{SW}(\Sigma(a_1,\cdots a_n),g_0) =2C_{a_1, \cdots ,a_n}, \;\;n=3,4
\]
where  $C_{a_1,\cdots, a_n}$ is the number of lattice points in the simplex
\[
\Delta(a_1, \cdots , a_n):=\left\{ (x_1, \cdots ,x_n)\in {\bR}^n\; ;\; x_i\geq 0,\;\sum_{i=1}^n\frac{x_i}{a_i}< \frac{1}{2}\left(n-2-\sum_{i=1}^n\frac{1}{a_i} \right)\right\}.
\]
 The arithmetic goal of this paper is the determination of  $C_{a_1,\cdots a_n}$ when $n=3,4$. This is a problem of independent interest since  the vertices of the simplex $\Delta(a_1, \cdots, a_n)$ are not lattice points and  none of the counting techniques using Riemann-Roch on toric manifolds seem to apply.    We will use instead a variation of a trick of Mordell; see \cite{Mo} or \cite{RG}.

 R. Fintushel and R. Stern  have shown in \cite{FS} that the Casson invariant of the Brieskorn  sphere $\Sigma(a,b,c)$  is  $\frac{1}{8}\sigma(a,b,c)$ where $\sigma(a,b,c)$ denotes the signature of the  Milnor  fiber of $\Sigma(a,b,c)$. This result was extended to arbitrary $\Sigma(a_1, \cdots , a_n)$ by Neuwmann-Wahl in \cite{NW}. The Kronheimer-Mrowka conjecture in this case    is equivalent to
\be
-2C_{a_1,\cdots , a_n} -\frac{1}{8}{\bf F}(g_0) = \frac{1}{8}\sigma(a_1,\cdots ,a_n).
\label{eq: km0}
\ee
The key objects in this paper will be the so called   {\em Dedekind-Rademacher} sums defined  for every coprime positive integers $h,k$ and any real numbers $x,y$ by
\[
s(h,k;x,y)=\sum_{\mu=0}^{k-1}\LL \frac{\mu+y}{k}\RR \LL \frac{ h(\mu +y) }{k}+x \RR
\]
where for any $r\in {\bR}$  we set
\[
((r))=\left\{
\begin{array}{rc}
0 & r\in {\bZ} \\
\{q\}-\frac{1}{2} & r\in {\bR}\setminus {\bZ}
\end{array}
\right.  \;\;(\{r\}:=r-[r]).
\]
Despite their apparent complexity these sums are computationally very friendly due mainly to the reciprocity law they satisfy (see \cite{Ra} or  the Appendix). The sums $s(h,k;0,0)$ are precisely the Dedekind sums $s(h,k)$ discussed in great detail in \cite{HZ} or \cite{RG}.

Our proof of (\ref{eq: km0}) is based on the following facts.

\noindent $\bullet$ According to Zagier (see \cite{HZ} or \cite{NW}) the signature $\sigma(a_1,\cdots ,a_n)$ can be expressed in terms of Dedekind sums.

\noindent $\bullet$ According to the computations in \cite{N} the quantity $\frac{1}{2}\eta_{dir}(g_0) +\frac{1}{8}\eta_{sign}(g_0)$ can be expressed in terms of Dedekind-Rademacher sums.

\noindent $\bullet$   For $n=3,4$ the number $C_{a_1,\cdots, a_n}$ can be expressed in terms of Dedekind-Rademacher sums.

As an arithmetic byproduct of the proof we obtain a   divisibility result for certain expressions involving Dedekind-Rademacher sums (see Remark \ref{rem: divi}).

The present paper  consists of three sections and an appendix.   In the first section we use the results of \cite{N} to express ${\bf F}$ in terms of Dedekind-Rademacher sums and to reduce the computation of $\chi_{SW}$ to a lattice point count. In the next section we  describe a variation of the  Mordell trick  which reduces the lattice point count to a certain arithmetic expression.  The third  section  describes this arithmetic expression in terms of Dedekind-Rademacher sums and completes the proof   of (\ref{eq: km}). For the reader's convenience we have included a brief appendix containing the basic  properties of   Dedekind-Rademacher sums used in this paper.

\medskip

\noindent {\bf Note} After this paper was completed we found out that the Kronheimer-Mrowka conjecture was proved independently and in its entire generality by Lim and Carey-Marcolli-Wang.  We believe the nice arithmetic behind the  Brieskorn homology spheres is interesting enough to  warrant a separate  treatment.

\tableofcontents

\section{Geometric preliminaries}

For pairwise coprime integers $a_1,\cdots ,a_n\geq 2$, $n\geq 3$ we denote by $\Sigma(\vec{a})$  the Brieskorn homology sphere $\Sigma(a_1,\cdots, a_n)$ with $n$ singular fibers (see \cite{JN} for a precise definition).
We orient $\Sigma(\vec{a})$ as the boundary of a complex manifold. $\Sigma(\vec{a})$ is a Seifert manifold. With respect to the above orientation it is a singular $S^1$ fibration over the  orbi-sphere $S^2(\vec{a})$  which has $n$ cone points of isotropies ${\bZ}_{a_i}$, $1\leq i \leq n$. This fibration has rational degree
\[
\ell =-\frac{1}{A}, \;\; A:=a_1a_2\cdots a_n.
\]
Set $b_i:=A/a_i$. The Seifert invariants $\vec{\beta}=(\beta_1,\cdots ,\beta_n)$  (normalized as in \cite{N}) are defined by
\[
\beta_ib_i\equiv -1 \,({\rm mod}\; a),\;\;0\leq \beta_i <a_i
\]
Set $q_i=\beta_i^{-1}=-b_i$ mod  $a_i$, $i=1,\ldots, n$. The canonical line bundle  of $S^2(\vec{a})$ has rational degree
\[
\kappa=(n-2)-\sum_{i=1}^n\frac{1}{a_i}.
\]
The universal  covering space of $\Sigma(\vec{a})$ is a Lie group $G=G(\vec{a})$ determined by
\[
G=\left\{
\begin{array}{lc}
SU(2) & \kappa  <0 \\
\tilde{PSL}_2({\bR}) & \kappa >0
\end{array}
\right.
\]
where $\tilde{PSL}_2({\bR})$ denotes the universal cover of $PSL_2({\bR})$.  Moreover $\Sigma(\vec{a})\cong \Gamma/G(\vec{a})$ where $\Gamma$ is a discrete subgroup of $G$. The natural left invariant metrics on $G$ (see \cite{S})  induce a natural metric $g_0$ on $\Sigma(\vec{a})$. All the  geometric quantities discussed  in the sequel are    computed with respect to  this metric and  for simplicity we will omit $g_0$ from  the various notations. Thus ${\bf F}(\vec{a})$ is ${\bf F}(g_0)$.

Set
\[
\rho:=\left\{\frac{\kappa}{2\ell}\right\}=\left\{
\begin{array}{cc}
\frac{1}{2} & A\; {\rm even}\\
0 & A \; {\rm odd}
\end{array}
\right.
\]
and define $\vec{\gamma}=(\gamma_1,\cdots,\gamma_n)$ by the equalities
\[
\gamma_i=m\beta_i\,({\rm mod}\; \alpha_i),\;\;1\leq i\leq n
\]
where $m$ is the integer
\[
m:=\frac{\kappa}{2\ell}-\rho =\frac{u-A-2\rho}{2},\;\;u:=\sum_ib_i.
\]
In\cite{N} we have proved the following. 

\medskip

\noindent $\bullet$ If  $A$ is even then
\be
{\bf F}(\vec{a}) =1-4\sum_i s(\beta_i, a_i) -4\sum_i\left(\, \LL\frac{q_i\gamma_i+\rho}{a_i}\RR+2s(\beta_i,a_i;\frac{\gamma_i+\beta_i\rho}{a_i},-\rho)\right)\label{eq: liviu1}
\ee
The above expression can be further  simplified using the identities
\be
s(\beta_i,a_i)=-s(b_i,a_i)
\label{eq: ele0}
\ee
\be
s(\beta_i,a_i;\frac{\gamma_i+\beta_i/2}{a_i}, -1/2)=-s(b_i,a_i,1/2,1/2)-\frac{1}{2}\LL\frac{q_i\gamma_i+1/2}{a_i}\RR
\label{eq: ele}
\ee
The identity (\ref{eq: ele0}) is elementary and can be safely left to the reader. The identity (\ref{eq: ele}) is proved in the Appendix.    Putting the above together we deduce that when $A$ is even we have
\be
{\bf F}(\vec{a})=1+4\sum_i s(b_i,a_i)+8\sum_is(b_i,a_i;1/2,1/2).
\label{eq: liv1}
\ee
$\bullet$ If $A$ is odd then
\be
{\bf F}(\vec{a})=1-\frac{1}{A}-4\sum_i s(\beta_i, a_i) -4\sum_{i=1}^n\left( 2s(\beta_i,a_i;\frac{\gamma_i+\beta_i\rho}{a_i},-\rho)+\LL\frac{q_i\gamma_i+\rho}{a_i}\RR\, \right).
\label{eq: liviu2}
\ee
Similarly, we have an identity
\be 
s(\beta_i, a_i;\frac{\gamma_i}{a_i},0)+\frac{1}{2}\LL\frac{q_i\gamma_i}{a_i}\RR=-s(b_i,a_i;1/2,1/2)
\label{eq: ele1}
\ee
and we deduce
\be
{\bf F}(\vec{a})=1-\frac{1}{A} +4\sum_i s(b_i,a_i) +8\sum_i s(b_i,a_i;1/2,1/2).
\label{eq: 2}
\ee

The signature $\sigma(\vec{a})$ of the Milnor  fiber of $\Sigma(\vec{a})$ can be expressed in terms of  Dedekind sums as well (see \cite{NW},  Sect.1)  and we have 
\be
\sigma(\vec{a})=-1-\frac{(n-2)A}{3}+\frac{1}{3A} +\frac{1}{3}\sum_i \frac{b_i}{a_i}-4\sum_i s(b_i,a_i).
\label{eq: hz}
\ee
We  deduce that
\be
{\bf F}(\vec{a})+\sigma(\vec{a})= -\frac{(n-2)A}{3}+\frac{\ve}{3A}+\frac{1}{3}\sum_i\frac{b_i}{a_i} +8\sum_is(b_i,a_i;1/2,1/2)
\label{eq: liv3}
\ee
where 
\[
{\ve}=\left\{
\begin{array}{rc}
1 & A\; {\rm even}
\\
-2 & A\;{\rm odd}
\end{array}
\right. .
\]

The description of $\chi_{SW}$ requires a bit more work. Introduce the simplex
\[
\Delta(\vec{a})=\{\vec{x}\in {\bZ}^n\; ;\; x_i\geq 0,\; \sum_i\frac{x_i}{a_i} < \kappa/2\}.
\]
For each $\vec{x}\in \Delta(\vec{a})$ set
\[
d(\vec{x})=\sum_i \left[\frac{x_i}{a_i}\right]
\]
and ${\cal S}_{\vec{x}}$= symmetric product of $d(\vec{x})$ copies of $S^2$.  Note that if $n=3,4$ then $d(\vec{x})=0$ for all $\vec{x}\in \Delta(\vec{a})$. 

The irreducible part  of the  adiabatic Seibert-Witten equations on $\Sigma(\vec{a})$ (studied in \cite{MOY}) and \cite{N1}) can be described as
\[
\modu_{\vec{a}}=\bigcup_{\vec{x}}\modu_{\vec{x}}
\]
where 
\[
\modu_{\vec{x}}=\modu^+_{\vec{x}}\cup \modu^-_{\vec{x}},\;\;\modu^\pm_{\vec{x}}\cong  {\cal S}_{\vec{x}}.
\]
Moreover, the  virtual dimensions of the  spaces of   finite energy gradient flows originating at the unique reducible solution and ending at one of the $\modu_{\vec{x}}$ are all odd.     Using the adiabatic argument in \S 3.3 of \cite{N} we deduce that   if all $d(\vec{x})$ are zero the Seiberg-Witten-Floer homology obtained using the usual Seiberg-Witten equations is isomorphic with the Seiberg-Witten-Floer homology obtained using the adiabatic equations.   Moreover,  all the even dimensional Betti numbers of  the Seiberg-Witten-Floer homology   are zero and  we deduce 
\be
\chi_{SW}(\vec{a}) =-2C_{\vec{a}}:=-2\# \Delta(\vec{a}).
\label{eq: count1}
\ee
In 1996 Kronheimer and Mrowka  conjectured that 
\[
\chi_{SW}(N,g)-\frac{1}{8}{\bf F}(N,g) =\lambda(N)
\]
for any  homology sphere $N$ and any generic metric $g$ on $N$.  Above, $\lambda(N)$  denotes the Casson invariant of  $N$.  It was shown in \cite{NW} that for the Brieskorn homology spheres $\Sigma(\vec{a})$
\be
\lambda(\Sigma(\vec{a}) )=\frac{1}{8}\sigma(\vec{a}).
\label{eq: nw}
\ee
The main result  of this paper is the following.

\begin{theorem}{\rm The Kronheimer-Mrowka conjecture is true for Brieskorn homology spheres $\Sigma(a_1, \cdots, a_n)$, $n=3,4$. According to   (\ref{eq: count1}), (\ref{eq: nw}) and(\ref{eq: liv3}) this is equivalent to}
\[
-16C_{\vec{a}}={\bf F}(\vec{a}) +\sigma(\vec{a})
\]
\be
=-\frac{(n-2)A}{3}+\frac{\ve}{3A}+\frac{1}{3}\sum_i\frac{b_i}{a_i} +8\sum_is(b_i,a_i;1/2,1/2).
\label{eq: km}
\ee
\label{th: liviu}
\end{theorem}

\begin{remark}{\rm As indicated in \cite{N}, Rohlin's theorem implies that the term ${\bf F}(\vec{a})$ is divisible by $8$.  The results of \cite{NW} show the  signature $\sigma(\vec{a})$ is also divisible by $8$. Thus the right-hand-side of (\ref{eq: liv3}) is an integer divisible by $8$. The  above theorem shows that ${\bf F}(\vec{a})+\sigma(\vec{a})$ is in effect divisible by 16!!!}
\label{rem: divi}
\end{remark}

\setcounter{equation}{0}

\section{The Mordell trick}
\setcounter{equation}{0}
Let $\vec{a}\in {\bZ}^n$ be as in the previous section. Denote by $\p=\p_{\vec{a}}$ the parallelipiped 
\[
\p := ([0,a_1-1]\times \cdots \times [0,a_n-1])\cap {\bZ}^n.
\]
When $n=3$ we will use the notation $\vec{a}=(a,b,c)$. Define $q:\p\ra {\bR}$ by
\[
q(\vec{x}) = \sum_i\frac{x_i+1/2}{a_i} =\sum_i\frac{x_i}{a_i}+\frac{u}{2A}.
\]
\begin{remark}{\rm (a) Suppose $n=3$, $\vec{a}=(a,b,c)$. Note that $q({\bf p})\in \frac{1}{2}{\bZ}$ for some ${\bf p}\in {\p}$ if and only if $abc$ is odd and 
\[
{\bf p}={\bf p}_0:=(\frac{a-1}{2}, \frac{b-1}{2}, \frac{c-1}{2}).
\]
In this case $q({\bf p}_0)=\frac{3}{2}$.

(b) Suppose $n=4$, $\vec{a}=(a_1,\cdots, a_4)$. Then $q({\bf p}\in {\bZ}$ for some ${\bf p} \in {\p}$ if and only if $A$ is odd and 
\[
{\bf p}={\bf p}_0=(\frac{a_1-1}{2},\cdots ,\frac{a_4-1}{2})
\]
in which case $q({\bf p}_0)=2$.}
\label{rem: parity}
\end{remark}
For every interval  $I\subset {\bR}$ we put
\[
N_I:=\# q^{-1}(I).
\]
Note that 
\[
C_{\vec{a}}=N_{(0,(n-2)/2)}.
\]
In particular, if $n=3$
\be
C_{a,b,c}=N_{(0,1/2)}.
\label{eq: mordell1}
\ee
while if $n=4$
\be
C_{a_1,\cdots, a_4}=N_{(0,1)}.
\label{eq: mord1}
\ee
For every $r\in {\bf R}$  define $\|r\|=[r+1/2]$ where $[\cdot ]$ denotes the integer part function. Note that $\|r\|$ is the integer closest to $r$.  We now discuss separately the two cases, $n=3$ and $n=4$

\medskip

\noindent $\bullet$ {\bf The case} $n=3$, $\vec{a}=(a,b,c)$. Immitating Mordell (see \cite{Mo} or \cite{RG}) we introduce the quantity
\[
\alpha := \sum_{\p}(\|q\|-1)(\|q\|-2).
\]
Observe that 
\be
\alpha =2N_{[0,1/2)}+2N_{[5/2,3)}=2N_{(0,1/2)}+2N_{(5/2,3)}.
\label{eq: mordell2}
\ee
The importance of the last equality follows from the following elementary result.
\begin{lemma}
\[
N_{(0,1/2)}=N_{(5/2,3)}
\]
\end{lemma}
\noindent {\bf Proof}\hspace{.3cm} Consider the involution
\[
\omega:\p \ra \p,\;\; (x,y,z)\mapsto (a-1-x,b-1-y,c-1-z).
\]
It has the property
\[
q(\omega({\bf p}))=3-q({\bf p})
\]
from which the lemma follows immediately.  $\Box$

\bigskip

Using the lemma and the equalities (\ref{eq: mordell1}), (\ref{eq: mordell2}) we deduce
\be
4C_{a,b,c}=\sum_{\p}(\|q\|-1)(\|q\|-2).
\label{eq: mordell3}
\ee
$\bullet$ {\bf The case} $n=4$, $\vec{a}=(a_1,\cdots,a_4)$. Arguing exactly as above we deduce
\be
4C_{\vec{a}}=N_{(0,1)}=\sum_{\bf p} ([q]-1)([q]-2).
\label{eq: mord3}
\ee
The proof of Theorem \ref{th: liviu} will be completed by providing an expression for the above sums in terms of Dedekind-Rademacher sums. This will be achieved in the next section following  the strategy of \cite{Mo} (see also \cite{RG}).

\section{The proof of Theorem \ref{th: liviu}}
\setcounter{equation}{0}

We will consider separately the two cases $n=3$ and $n=4$.

\subsection{The case $n=3$}
Set $\vec{a}=(a,b,c)$ so that $A=abc$, $u=ab+bc+ca$. We will distinguish two cases: $A$ is even and $A$ is odd.

\medskip

\noindent $\bullet$ $A$ {\bf is even} In this case $q({\bf p}) + \frac{1}{2}\not\in {\bZ}$ so that
\[
\|q\|=q+1/2 -\{q+1/2\} =q-((q+1/2)).
\]
The sum  can be rewritten as
\[
\sum_{\p}(q-((q+1/2))-1)(q-((q+1/2))-2)
\]
\[
=\sum_{\p}q^2 -3\sum_{\p}q+2\sum_{\p}1 
\]
\[
-2\sum_{\p}q((q+1/2)) +\sum_{\p}((q+1/2))^2+3\sum_{\p}((q+1/2)).
\]
We compute each of these 6 sums separately.
\[
\sum_{\p} 1=\# \p =abc.
\]
\[
\sum_{\p}q=\sum_{\p}\frac{x+1/2}{a}+\sum_{\p}\frac{y+1/2}{b}+\sum_{\p}\frac{z+1/2}{c}
\]
\[
=\frac{bc}{2a}\sum_{x=0}^{a-1}(2x+1)+\frac{ca}{2b}\sum_{y=0}^{b-1}(2y+1)+\frac{ab}{2c}\sum_{z=0}^{c-1}(2z+1)
\]
\[
=\frac{3abc}{2}.
\]
\[
\sum_{\p}q^2= \sum_{\p}\left(\frac{x+1/2}{a}\right)^2\sum_{\p}\left(\frac{y+1/2}{b}\right)^2 +\sum_{\p}\left(\frac{z+1/2}{c}\right)^2
\]
\[
+ 2c \left( \sum_{x=0}^{a-1}\frac{x+1/2}{a} \right)  \left(\sum_{y=0}^{b-1}\frac{y+1/2}{b}\right) + 2b \left( \sum_{x=0}^{a-1}\frac{x+1/2}{a} \right)  \left(\sum_{z=0}^{c-1}\frac{z+1/2}{c}\right) 
\]
\[
 +2a \left(\sum_{z=0}^{c-1}\frac{z+1/2}{c}\right)  \left(\sum_{y=0}^{b-1}\frac{y+1/2}{b}\right) .
\]
Using  basic properties of Bernoulli polynomials (see \cite{Ra1}) we deduce
\[
\sum_{x=0}^{a-1}\left(\frac{x+1/2}{a}\right)^2 =\frac{1}{3a^2}\left(B_3(a+\frac{1}{2})-B_3(\frac{1}{2})\right)
\]
where
\[
B_3(t)=\frac{t(2t-1)(t-1)}{2}
\]
is the third Bernoulli polynomial. Note that $B_3(1/2)=0$ and 
\[
B_3(t+1/2)=t(t^2-1/4).
\]
Using the identity
\[
\frac{1}{n}\sum_{k=0}^{n-1}\frac{k+1/2}{n}=\frac{n}{2}
\]
we conclude
\[
\sum_{\p}q^2=\frac{5abc}{2}-\frac{1}{12}\left( \frac{bc}{a}+\frac{ca}{b}+\frac{ab}{c}\right).
\]
Next
\[
\sum_{\p}((q+1/2))=\sum_{\p} \LL \frac{x}{a}+\frac{y}{b}+\frac{z}{c}+\frac{u}{2abc}\RR
\]
\[
=\sum_{k=0}^{abc-1}\LL \frac{k}{abc}+\frac{u}{2abc}\RR
\]
According to the Kubert identity (\ref{eq: kubert}) in the Appendix  the last sum is equal to $((u/2))$ which is zero. Thus
\[
\sum_{\p}((q+1/2))=0.
\]
The sum $\sum q((q+1/2))$ requires a bit more work. Note first that it decomposes as
\[
\sum_{x=0}^{a-1}\frac{x+1/2}{a} \sum_{y,z}\LL\frac{x}{a}+\frac{y}{b}+\frac{z}{c}+\frac{u+abc}{2abc}\RR
\]
\[
+\sum_{y=0}^{b-1}\frac{y+1/2}{b} \sum_{z,x}\LL\frac{x}{a}+\frac{y}{b}+\frac{z}{c}+\frac{u+abc}{2abc}\RR
\]
\[
+ \sum_{z=0}^{c-1}\frac{z+1/2}{c}\sum_{x,y}\LL\frac{x}{a}+\frac{y}{b}+\frac{z}{c}+\frac{u+abc}{2abc}\RR
\]
\[
=S_1+S_2+S_3.
\]
We    describe in detail the computation of $S_1$. The other two sums  are  entirely similar. Note first that
\[
\sum_{y,z}\LL\frac{x}{a}+\frac{y}{b}+\frac{z}{c}+\frac{u+abc}{2abc}\RR=\sum_{y,z}\LL\frac{y}{b}+\frac{z}{c}+\frac{x}{a}+\frac{u+abc}{2abc}\RR
\]
\[
=\sum_{k=0}^{bc-1}\LL \frac{k}{bc}+\frac{x}{a}+\frac{u+abc}{2abc}\RR
\]
(use the Kubert identity (\ref{eq: kubert}))
\[
=\LL \frac{bcx}{a}+\frac{u+abc}{2a}\RR.
\]
\[
=\LL \frac{bc(x+1/2)}{a} +\frac{bc+b+c}{2}\RR =\LL \frac{bc(x+1/2)}{a}+\frac{1}{2}\RR.
\]
We  conclude
\[
S_1=\sum_{x=0}^{a-1}\frac{x+1/2}{a}\LL \frac{bc(x+1/2)}{a}+\frac{1}{2}\RR
\]
\[
=\sum_{x=0}^{a-1}\LL\frac{x+1/2}{a}\RR  \LL\frac{bc(x+1/2)}{a}+\frac{1}{2}\RR -\frac{1}{2}\sum_{x=0}^{a-1}\LL\frac{bc(x+1/2)}{a}+\frac{1}{2}\RR
\]
(use the Kubert identity)
\[
=\sum_{x=0}^{a-1}\LL\frac{x+1/2}{a}\RR  \LL\frac{bc(x+1/2)}{a}+\frac{1}{2}\RR 
\]
\[
=s(bc,a;1/2,1/2).
\]
Hence
\[
\sum_{\p}q((q+1/2))=s(bc,a;1/2,1/2)+s(ca,b;1/2,1/2)+s(ab,c;1/2,1/2).
\]
Finally
\[
\sum_{\p}((q+1/2))^2=\sum_{\p} \LL\frac{x}{a}+\frac{y}{b} +\frac{z}{c} +\frac{u+abc}{2abc}\RR^2
\]
\[
=\sum_{k=0}^{abc-1}\LL \frac{k+\frac{u+abc}{2}}{abc}\RR^2
\]
( use the fact that $u+abc$ is odd in this case)
\[
=\sum_{k=0}^{abc-1}\LL\frac{k+1/2}{abc}\RR^2 =s(1,abc;0,1/2)
\]
(use  (\ref{eq: quad})
\[
= \frac{abc}{12}-\frac{1}{12 abc}.
\]
Putting together  the above information  we deduce that if $abc$ is even then
\[
4C_{a,b,c}=\frac{abc}{12}-\frac{1}{12abc}-\frac{1}{12}\left(\frac{bc}{a}+\frac{ca}{b}+\frac{ab}{c}\right)
\]
\be
-2(s(bc,a;1/2,1/2)+s(ca,b;1/2,1/2)+s(ab,c;1/2,1/2)).
\label{eq: mordell4}
\ee
The identity (\ref{eq: km}) is now obvious.

\medskip

\noindent $\bullet$ $A$ {\bf is odd} In this case, using Remark \ref{rem: parity} we deduce 
\[
\|q(p)\|=\left\{
\begin{array}{lr}
q-((q+1/2)) & {\bf p}\neq {\bf p}_0 \\
q-(q+1/2))+\frac{1}{2} & {\bf p}={\bf p}_0
\end{array}
\right. .
\]
Thus
\[
(\|q\|-1)(\|q\|-2)=\left\{
\begin{array}{lr}
(q-((q+1/2))-1)(q-((q+1/2))-2) &  {\bf p}\neq {\bf p}_0 \\
(q-((q+1/2))-1/2)(q-((q+1/2))-3/2)  & {\bf p}={\bf p}_0
\end{array}
\right.
\]
Hence
\[
4C_{a,b,c}= \sum_{\p}(q-((q+1/2))-1)(q-((q+1/2))-2)
\]
\[
+(q-((q+\frac{1}{2}))-\frac{1}{2})(q-((q+\frac{1}{2}))-\frac{3}{2})\!\mid_{{\bf p}_0}-(q-((q+\frac{1}{2}))-1)(q-((q+\frac{1}{2}))-2)\!\mid_{{\bf p}_0}
\]
\be 
=\sum_{\p}(q-((q+1/2))-1)(q-((q+1/2))-2) +\frac{1}{4}.
\label{eq: mordell5}
\ee
The above sum can be computed exactly as in the even case with one notable difference namely
\[
\sum((q+1/2))^2 =\sum_{k=0}^{abc-1}\LL \frac{k +\frac{u+abc}{2}}{abc}\RR
\]
($u+abc$ is even)
\[
=\sum_{k=0}^{abc-1}\LL\frac{k}{abc}\RR = s(1,abc;0,0) =\frac{abc}{12}+\frac{1}{6abc}-\frac{1}{4}.
\]
Thus, when $abc$ is odd we have
\[
4C_{a,b,c}=\frac{abc}{12}+\frac{1}{6abc}-\frac{1}{12}\left(\frac{bc}{a}+\frac{ca}{b} +\frac{ab}{c}\right)
\]
\[
-2(s(bc,a;1/2,1/2)+s(ca,b;1/2,1/2)+s(ab,c;1/2,1/2)).
\]
This completes the proof of Theorem \ref{th: liviu} when $n=3$

\subsection{The case $n=4$}
We follow a similar strategy with some obvious modifications. Set $\vec{a}=(a_1,\cdots, a_4)$, $A=4$, $u=b_1+\cdots +b_4$ and
\[
S_{\vec{a}}= \sum_{ {\p}_{ \vec{a} } }([q]-1)([q]-2).
\]
As in the previous subsection  will distinguish two situations.

\medskip 

\noindent $\bullet$ $A$ {\bf is even} Note that  for every ${\bf p}\in {\p}$ we have $q({\bf p}) \not\in {\bZ}$ so that
\[
[q]=q-((q))-1/2.
\]
Thus
\[
S_{\vec{a}}=\sum_{\p} (q-((q))-3/2)(q-(q))-5/2)
\]
\[
=\sum_{\P}(q^2-4q+15/4)-2\sum_{\p}q((q))+\sum_{\p}((q))^2+4\sum_{\p}((q)).
\]
The computation of the above terms follows the same pattern  as in the previous subsection.
\[
\sum_{\p}((q))=0.
\]
\[
\sum_{\p}15/4 =15\#{\p}/4=15A/4.
\]
\[
\sum_{\p}q= \sum_{i=1}^4b_i\sum_{x_i=0}^{a_i-1}\frac{x_i+1/2}{a_i}=\sum_{i=1}^4b_ia_i=2A.
\]
\[
\sum_{p}q^2=\sum_{i=1}^4b_i\sum_{x_i=0}^{a_i-1}\left(\frac{x_i+1/2}{a_i}\right)^2 
\]
\[
+2 \sum_{i<j}\frac{A}{a_ia_j}\left(\sum_{x_i=0}^{a_i-1}\frac{x_i+1/2}{a_i}\right)\left(\sum_{x_j=0}^{a_j-1}\frac{x_j+1/2}{a_j}\right)
\]
\[
=\sum_{i=1}^4\frac{b_ia_i(a_i^2-1/4)}{3a_i^2}+\sum_{1\leq i<j\leq 4}\frac{A}{2}
\]
\[
=\sum_{i=1}^{4}\left(\frac{A}{3}-\frac{b_i}{12a_i}\right) +3A
=\frac{13A}{3}-\frac{1}{12}\sum_i \frac{b_i}{a_i}.
\]
When $A$ is even $u$ is odd  and  we have
\[
\sum_{\p}((q))^2 =\sum_{k=0}^{A-1}\LL \frac{k+u/2}{A}\RR^2
=s(1,A;0,1/2)=\frac{A}{12}+\frac{1}{2A}.
\]
Finally,
\[
\sum_{\p} q((q)) = S_1+\cdots +S_4
\]
where 
\[
S_1=\sum_{x_1=0}^{a_1-1}\frac{x_1+1/2}{a_1}\cdot \sum_{x_2,x_3,x_4}\LL\frac{x_2}{a_2}+\frac{x_3}{a_3}+\frac{x_4}{a_4}+\frac{x_1}{a_1}+\frac{u}{2A}\RR.
\]
$S_2,S_3, S_4$ are defined similarly. To compute $S_1$ note that
\[
\sum_{x_2,x_3,x_4}\LL\frac{x_2}{a_2}+\frac{x_3}{a_3}+\frac{x_4}{a_4}+\frac{x_1}{a_1}+\frac{u}{2A}\RR
\]
\[
=\sum_{k=0}^{b_1-1}\LL\frac{k}{b_1}+\frac{x_1}{a_1}+\frac{u}{2A}\RR
\]
(use the Kubert identity)
\[
=\LL\frac{b_1x_1}{a_1}+\frac{u}{2A}\RR =\LL \frac{b_1(x_1+1/2)}{a_1}+\frac{u-b_1}{2a_1}\RR
\]
($(u-b_1)/a_1$ is odd)
\[
= \LL\frac{b_1(x_1+1/2)}{a_1}+ \frac{1}{2}\RR.
\]
Thus
\[
S_1=\sum_{x_1}{a_1-1}\frac{x_1+1/2}{a_1} \LL\frac{b_1(x_1+1/2)}{a_1}+ \frac{1}{2}\RR
\]
and we deduce as in the previous subsection that
\[
S_1= s(b_1,a_1;1/2,1/2).
\]
By adding all the above together we deduce that if $A$ is odd then
\[
4C_{\vec{a}}=S_{\vec{a}}=\frac{A}{6}-\frac{1}{12A} -\frac{1}{12}\sum_i\frac{b_i}{a_i} -8\sum_is(b_i,a_1;1/2,1/2).
\]
The identity (\ref{eq: km}) is now obvious.

\medskip

\noindent $\bullet$ $A$ {\bf is odd} In this case $u$ is even. Arguing as in the previous subsection we deduce
\[
S_{\vec{a}}=\sum_{\p}(q-((q))-3/2)(q-((q))-5/2) +\frac{1}{4}.
\]
The only term in the previous computations which is influenced  by the parity of $A$ is
\[
\sum_{\p}((q))^2 =\sum_{k=0}^{A-1}\LL \frac{k+u/2}{A}\RR^2
=s(1,A)
\]
\[
= \frac{A}{12}+\frac{1}{6A}-\frac{1}{4}.
\]
Putting together all the terms  we obtain again the identity (\ref{eq: km}).  The Theorem \ref{th: liviu} is proved. $\Box$

\appendix 
\section{Basic facts concerning Dedekind-Rademacher sums}
\setcounter{equation}{0}
In \cite{Ra} Rademacher consider for every pair of coprime  integers $h,k$ and any real numbers  $x,y$  the following generalization of the classical Dedekind sums
\[
s(h,k;x,y)=\sum_{\mu=0}^{k-1}\LL \frac{\mu+y}{k}\RR \LL \frac{ h(\mu +y) }{k}+x \RR.
\]
A simple computations shows that $s(h,k; x,y)$ depends only on $x,y$ mod 1.   When $h=1$ and $x=0$  one can prove (see \cite{Ra})
\be
s(1,k;0,y)=\left\{
\begin{array}{lr}
\frac{k}{12}+\frac{1}{6k}-\frac{1}{4} & y \in {\bZ} \\
                                     &   \\
\frac{k}{12} +\frac{1}{k}B_2(\{y\})  & y \in {\bR}\setminus {\bZ}
\end{array}
\right.
\label{eq: quad}
\ee
where $B_2(t)=t^2-t+1/6$ is the second Bernoulli polynomial. 

Perhaps the  most important property of these Dedekind-Rademacher sums is their reciprocity law which makes them    computationally very friendly. To formulate it we must distinguish two cases.

\noindent  $\bullet$  Both $x$ and $y$ are integers.  Then
\be
 s(\beta, \alpha; x,y) +s(\alpha, \beta; y,x) =-\frac{1}{4}+\frac{\alpha^2+\beta^2+1}{12\alpha \beta}.
\label{eq: rec1}
\ee
$\bullet$  $x$ and/or $y$ is not an integer. Then
\be
s(\beta,\alpha; x, y)+s(\alpha, \beta; y,x)=((x))\cdot ((y)) + \frac{\beta^2\psi_2(y) +\psi_2(\beta y+\alpha x) +\alpha^2\psi_2(x)}{2\alpha \beta} 
\label{eq: rec2}
\ee
where $\psi_2(x):= B_2(\{x\})$. 

An important ingredient behind  the reciprocity law is the following identity (Lemma 1 in \cite{Ra})
\be
\sum_{\mu=0}^{k-1}\LL\frac{\mu+w}{k}\RR =((w))\;\;\forall w\in {\bR}.
\label{eq: kubert}
\ee
Following the terminology in \cite{Mi} we will call the above equality the {\em  Kubert  identity}.

We conclude with a proof of the identity (\ref{eq: ele}). For simplicity we consider only the case $n=3$ and $i=1$. Set $\vec{a}=(a,b,c)$. Thus $A=abc$ is even, $u=bc+ca+ab$ and $b_1=bc$.  For arbitrary $n$ the proof is only notationally more complicated.

The proof  of (\ref{eq: ele}) goes as follows.
\[
s(\beta_1,a;\frac{\gamma_1+\beta_1/2}{a},-1/2)  = \sum_{x=0}^{a-1} \LL\frac{x-1/2}{a}\RR \LL\frac{\beta_1x +\gamma_1}{\alpha_1}\RR
\]
($\gamma_1=\beta_1(u-abc-1)/2$ mod $a$)
\[
= \sum_{x=0}^{a-1} \LL\frac{x-1/2}{a}\RR \LL \frac{\beta_1(x-\frac{abc-u+1}{2}}{a}\RR
\]
($y:=x-\frac{abc-u+1}{2}$ mod $a$)
\[
=\sum_{y=0}^{a-1} \LL \frac{y +\frac{abc-u+1}{2} -1/2}{a}\RR \LL \frac{\beta_1 y}{a}\RR 
\]
( use $ z=-bc y$ mod $a$  and $\beta_1 bc \equiv 1$ mod $a_1$)
\[
=-\sum_{z=0}^{a-1}\LL \frac{bcz-\frac{abc-u}{2}}{a}\RR \LL \frac{z}{a}\RR
\]
\[
=-\sum_{z=0}^{a-1} \LL \frac{bc(z+1/2)}{a}+\frac{b+c-bc}{2}\RR \LL \frac{z}{a}\RR
\]
\[
=- \sum_{z=0}^{a-1}\LL \frac{bc(z+1/2)}{a}+\frac{1}{2} \RR \LL \frac{z}{a}\RR.
\]
At this point we use the elementary identity
\[
\LL \frac{z}{a}\RR =\LL \frac{z+1/2}{a}\RR -\frac{1}{2a} +\frac{1}{2}\delta (z)
\]
where
\[
\delta(z)=\left\{
\begin{array}{ll}
1 & z\equiv 0 \,({\rm mod}) a \\
0 & {\rm otherwise}
\end{array}
\right.
\]
We deduce
\[
s(\beta_1,\alpha_1;\frac{\gamma_1+\beta_1/2}{\alpha_1},-1/2)  = - \sum_{z=0}^{a-1}\LL \frac{bc(z+1/2)}{a}+\frac{1}{2} \RR \LL \frac{z+1/2}{a}\RR
\]
\[
+\frac{1}{2a}\sum_{z=0}^{a-1} \LL \frac{bc(z+1/2)}{a}+\frac{1}{2}\RR  -\frac{1}{2}\LL\frac{bc}{2a}+\frac{1}{2}\RR.
\]
The Kubert identity shows that the second sum above vanishes. Also
\[
\LL \frac{q_1\gamma_1 +1/2}{\alpha_1}\RR =\LL \frac{\frac{u-abc}{2}}{a}\RR =\LL \frac{b+c-bc}{2}+\frac{bc}{2a}\RR
\]
\[
=\LL
\frac{bc}{2a}+\frac{1}{2}\RR.
\]
The identity (\ref{eq: ele}) is proved.  The proof  of (\ref{eq: ele1})   is similar and is left to the reader.

\end{document}